\documentclass[a4paper,12pt]{amsart}
\usepackage{amsmath,amssymb,amsthm,ascmac,fancybox}   
\usepackage{bm}                
\usepackage{tabularx}          
\usepackage[dvipdfm]{graphicx} 
\allowdisplaybreaks[3]
\newcommand{\sa}   [1]{\mathcal{#1}}                       
  
\newcommand{\norm} [2]{\left\| #1 \right\|_{#2}}               
\newcommand{\ab}  [1]{\left| #1\right|} 				
\newcommand{\qw}   [0]{\par{\noindent}}                    

\newcommand{\N}{\mathbb{N}}			
\newcommand{\R}{\mathbb{R}}			

\newcommand{\alm}{\mathrm{a.e.}\,}		
\newcommand{\kko}  [1]{\left( #1\right)}  
\newcommand{\kkko}  [1]{\left\{ #1\right\}}	

\newcommand{\ds}{\displaystyle}

\title{An estimate for derivative \\ of the de la Vall\'{e}e Poussin mean}
\author{Kentaro\ Itoh, Ryozi\ Sakai and Noriaki\ Suzuki}	

\begin{document}
\
\begin{abstract}
The de la Vall\'{e}e Poussin mean for exponential weights on $(-\infty, \infty)$ was 
investigated in [6]. In the present paper we discuss its derivatives. 
An estimate for the Christoffel function plays an important role.  
\end{abstract}

\maketitle

\section{{Introduction}}
Let $\R=(-\infty,\infty)$. We consider an exponential weight
\[ w(x)=\exp(-Q(x)) \]
on $\R$, where $Q$ is an even and nonnegative function on $\R$.\ Throughout this
paper we always assume that $w$ belongs to a relevant class $\sa{F}(C^2+)$ (see section 2). A function $T = T_{w}$ defined by
\[ T(x):=\frac{xQ'(x)}{Q(x)},~~x\neq 0 \leqno{(1.1)} \]
is very important. We call $w$ a Freud-type weight if $T$ is bounded, and otherwise, $w$ is
called an Erd{\"{o}}s-type weight. For $x >0$, the Mhaskar-Rakhmanov-Saff number (MRS number) $a_x = a_{x}(w)$ of $w=\exp(-Q)$ is defined by a positive root of the equation
\[ x=\frac{2}{\pi}\int_0^1\frac{a_xuQ^\prime(a_xu)}{(1-u^2)^{1/2}}du. \leqno{(1.2)} \]
When $w = \exp(-Q) \in \sa{F}(C^2+)$, $Q'$ is positive and increasing on $(0, \infty)$, so that
$$
\lim_{x \to \infty} a_x = \infty  \ \ \mbox{and} \ \ 
\lim_{x \to +0 }a_{x} = 0 \leqno(1.3)
$$
and 
$$
\lim_{x \to \infty} \frac{a_x}{x} = 0 \ \ \mbox{and} \ \ \lim_{x \to +0} \frac{a_x}{x} = \infty \leqno(1.4)
$$ hold. Note that those convergences are all monotonically.

Let $\{p_n\}$ be orthogonal polynomials for a weight $w$, that is, $p_n$ is the polynomial of degree $n$ such that
$$ \int_\R p_n(x)p_m(x)w(x)^{2}dx = \delta_{mn}. $$
Note that when $w(x) = \exp(-|x|^{2})$, then $\{p_{n}\}$ are Hermite polynomials. 

For $1\le p\le\infty$, we denote by $L^p(I)$ the usual $L^p$ space on an interval $I$ in  $\R$.  For a function $f$ with $fw\in L^p(\R)$, we set
\[ s_n(f)(x):=\sum_{k=0}^{n-1}b_k(f)p_{k}(x) ~~\mbox{where}~~b_k(f)=\int_\R f(t)p_k(t)w(t)^{2}dt \]
for $n\in\N$ (the partial sum of Fourier series). The de la Vall\'{e}e Poussin mean $v_n(f)$ of $f$ is defined by
\[ v_n(f)(x):=\frac{1}{n}\sum_{j=n+1}^{2n}s_j(f)(x). \]

In [6], we proved the following; Let $1\le p\le\infty$ and $w\in\sa{F}(C^2+)$. Assume that $T(a_n)\le Cn^{2/3-\delta}$ 
for some $0<\delta\le 2/3$ and $C>1$. Then there exists another constant $C>1$ such that if $fw \in L^{p}(\R)$, 
then
\[ \|v_n(f) \frac{w}{T^{1/4}}\|_{L^p(\R)}\le C\|fw\|_{L^p(\R)} \leqno{(1.5)} \]
holds for all $n \in \N$, and if  $T^{1/4}fw\in L^p(\R)$, then 
$$ \|v_n(f)w\|_{L^p(\R)}\le C\|T^{1/4}fw\|_{L^p(\R)} \leqno{(1.6)}$$
holds  for all $n \in \N$. 
It is also known that 
\[ \left\|P^\prime\frac{w}{T^{1/2}}\right\|_{L^p(\R)}\le C\kko{\frac{n}{a_n}}\|Pw\|_{L^p(\R)}, \leqno{(1.7)} \]
for all $P\in\sa{P}_n$, where $\sa{P}_n$ is the set of all polynomials of degree at most $n$ (see [5, Theorem 6.1]). 
Since $v_{n}(f) \in \sa{P}_{2n}$, combining  (1.6) with (1.7), we have  
$$
\left\|v_n^\prime(f)\frac{w}{T^{1/2}}\right\|_{L^p(\R)} \leq C \kko{\frac{n}{a_{n}}} \|T^{1/4}fw\|_{L^p(\R)} \leqno(1.8)
$$
with some $C >1$. Here we use the fact that $a_{n}$ and $a_{2n}$ are comparable (see Lemma 2.1 (1) below).
The inequality (1.8) suggests us the following: 
if $fw \in L^{p}(\R)$, then
\[ \norm{v_n^\prime(f) \frac{w}{T^{3/4}}}{L^p(\R)} \le C\kko{\frac{n}{a_n}}\|fw\|_{L^p(\R)}\leqno{(1.9)} \]
and, 
if $T^{3/4}fw\in L^p(\R)$, then
\[ \norm{v_n^\prime(f)w}{L^p(\R)} \le C\kko{\frac{n}{a_n}}\|T^{3/4}fw\|_{L^p(\R)}\leqno{(1.10)} \]
holds?

In the present paper, we will show that (1.9) holds for all $1 \leq p \leq \infty$ and (1.10) is true for $2 \leq p \leq \infty$ at the least. More generally, 
as for the $j$th derivative $v_{n}^{(j)}(f)$ of $v_n(f)$,  the following theorems are established. 

\vspace{2ex}

{\bf Theorem 1.1.} Let $k \geq 2$ be an integer and let $w\in\sa{F}_\lambda(C^{4}+)$ with $0<\lambda<(k+3)/(k+2)$, and let $1 \leq p \leq \infty$.  
Then there exists a constant $C >1$ such that if $1 \leq j \leq k$, and if $fw\in L^p(\R)$, then  
\[ \|v_n^{{(j)}}(f) \frac{w}{T^{(2j+1)/4}}\|_{L^p(\R)} \le C\kko{\frac{n}{a_n}}^j\|fw\|_{L^p(\R)}  \leqno{(1.11)} \]
holds for all $n \in \mathbb{N}$. 

\vspace{2ex}
The definition of a class $\sa{F}_{\lambda}(C^{4}+)$ is given in section 2.
\vspace{2ex}

{\bf Theorem 1.2.} Let $ k$ and $w$ be as in Theorem 1.1, and let $2 \leq p \leq \infty$.  
Then there exists a constant $C >1$ such that if $1 \leq j \leq k$, and if $T^{(2j+1)/4}fw\in L^p(\R)$, then  
\[ \|v_n^{{(j)}}(f) w\|_{L^p(\R)} \le C\kko{\frac{n}{a_n}}^j\|T^{(2j+1)/4}fw\|_{L^p(\R)}  \leqno{(1.12)} \]
holds for all $n \in \mathbb{N}$. 

\vspace{2ex}

{\bf Theorem 1.3.} Let $ k$ and $w$ be as in Theorem 1.1, and let $1 \leq p \leq 2$.  Then there exists a constant $C >1$ such that for every $1 \leq j \leq k$ and every  $T^{(2j+1)/4}fw\in L^2(\R)$, we have 

\[ \|v_n^{(j)}(f)w\|_{L^p(\R)} \le C\kko{\frac{n}{a_n}}^ja_n^{(2-p)/2p}\|T^{(2j+1)/4}fw\|_{L^2(\R)}  \leqno{(1.13)} \]
for all $n \in \mathbb{N}$. 

\vspace{2ex}

We note that when $w$ is a Freud-type weight, then $1 \leq T \leq C$, so that, 
$$
 \|v_n^{{(j)}}(f) w\|_{L^p(\R)} \le C\kko{\frac{n}{a_n}}^j\|fw\|_{L^p(\R)} \leqno (1.14)
 $$ follows from Theorem 1.1. In [3, Chapter 3], Mhasker discussed the 1st derivative of the de la Vall\'{e}e Poussin mean for Freud-type weights. Our contribution is to deal with not only Freud-type but also Erd\"{o}s-type weights. 
In the proofs of above theorems, we use Mhasker's argument. In addition, there are two keys: one is to use mollification of exponential weights (see Lemma 2.4 below) which was obtained in [5], and another
is to estimate the Christoffel functions which are done in section 3. Unfortunately, we do not know whether (1.12) holds true or not for $1 \leq p <2$, however, we will give another estimate which holds for all $1 \leq p \leq \infty$ in section 4.
A related inequality  to (1.13) is also given in section 6.

Throughout this paper, we write 
$f(x) \sim  g(x)$ for a subset $I \subset \R$ if there exists a constant $C\geq1$ such that $f(x)/C \leq g(x) \leq Cf(x)$ holds for all $x \in I$. 
Similarly, $a_{n} \sim b_{n}$ means that $a_{n}/C \leq b_{n} \leq Ca_{n}$ holds for all $n \in \N$. We will use 
the same letter $C$ to denote various positive constants; it may vary even within a line. Roughly speaking, $C >1$ implies  that $C$ is sufficiently large, and differently, $C >0$ means $C$ is a sufficiently small positive number. 

\section{{Definitions and Lemmas}}

We say that an exponential weight $w=\exp(-Q)$ belongs to class $\sa{F}(C^2+)$, when $Q:\R\to[0,\infty)$ is a continuous and even function and satisfies the following conditions:\\
$~~$(a)$~~Q^\prime(x)$ is continuous in $\R$ and  $Q(0)=0$.\\
$~~$(b)$~~Q^{\prime\prime}(x)$ exists and is positive in $\R\setminus\{0\}$.\\
$~~$(c)$\ds~~\lim_{x\to\infty}Q(x)=\infty.$\\
$~~$(d)$~~$The function $T$  in (1.1) is quasi-increasing in $(0,\infty)$(i.e. there exists $C>1$ such that $T(x) \le C T(y)$ whenever $0<x<y$), and there exists $\Lambda\in\R$ such that
\begin{align*} 
T(x)\ge \Lambda>1,~~x\in \R\setminus\{0\}.
\end{align*}
$~~$(e)$~~$There exists $C>1$ such that
$$
\frac{Q''(x)}{|Q'(x)|}\le C\frac{|Q'(x)|}{Q(x)},~~\alm~x\in \R.
$$
Moreover, if there also exist a compact subinterval $J (\ni 0)$ of $\R$, and $C>1$ such that
\begin{align*} 
 C \frac{Q''(x)}{|Q'(x)|}\ge \frac{|Q'(x)|}{Q(x)},~~\alm~x\in \R\setminus J.
\end{align*}

\vspace{2ex}

Let $\lambda >0$. We write $w \in \sa{F}_\lambda(C^2 +)$ if  there exist $K>1$ and $C>1$ such that for all $|x| \geq K$, 
$$ \frac{|Q^\prime(x)|}{Q(x)^\lambda}\le C \leqno(2.1)
$$
holds.  We also  write $w \in \sa{F}_\lambda(C^3 +)$, if $Q \in C^3(\R \setminus \{0\})$ and 
\[ \ab{\frac{Q^{(3)}(x)}{Q''(x)}}\le C \ab{\frac{Q^{\prime\prime}(x)}{Q^{\prime}(x)}} \ \ 
\mbox{and} \ \  \frac{|Q^\prime(x)|}{Q(x)^\lambda}\le C \]
holds for every $|x| \geq K$.   Moreover, we write $w \in \sa{F}_\lambda(C^4 +)$,  if $Q \in C^4(\R \setminus \{0\})$ and 
\[     \ab{\frac{Q^{(3)}(x)}{Q''(x)}} \sim   \ab{\frac{Q^{''}(x)}{Q'(x)}}, \ \  \ab{\frac{Q^{(4)}(x)}{Q^{(3)} (x)}}\le C \ab{\frac{Q^{\prime\prime}(x)}{Q^{\prime}(x)}} \ \ 
\mbox{and} \ \  \frac{|Q^\prime(x)|}{Q(x)^\lambda}\le C \]
holds for every $|x| \geq K$.  Clearly $\sa{F}_\lambda (C^4 +) \subset \sa{F}_\lambda (C^3 +) \subset \sa{F}_\lambda (C^2 +) \subset \sa{F}(C^2 +)$.

A typical example of Freud-type weight is $w=\exp(-|x|^\alpha)$ with $\alpha>1$. 
It belongs to $\sa{F}_\lambda(C^{4}+)$ for $\lambda =1$. 
For 
$u\ge 0,\ \alpha>0$ with $\alpha+u>1$ and $l\in\N$, we set 
\[ Q(x):=|x|^u(\exp_l(|x|^\alpha)-\exp_l(0)), \]
where $\exp_l(x) :=\exp(\exp(\exp(\cdots(\exp x))))$ ($l$-times). Then $w=\exp(-Q(x))$
is an Erd{\"{o}}s-type weight, which  belongs to $\sa{F}_\lambda(C^{4}+)$ for $\lambda >1$ (see[1]).

\vspace{2ex}

In the following lemmas we fix $w \in \sa{F}(C^2 +)$. 

\vspace{2ex}

{\bf Lemma 2.1.}  Fix $L>0$. Then we have 

(1)  $a_t \sim  a_{Lt}$ on $t>0$ (see [2, Lemma 3.5 (a)]). 

(2) $Q(a_t) \sim  Q(a_{Lt})$, $Q^\prime(a_t)\sim Q^\prime(a_{Lt})$ and $T(a_{Lt})\sim T(a_t)$ 
on $t>0$ (see [2, Lemma 3.5 (b)]). 

(3)  $\ds\frac{1}{T(a_t)}\sim \ab{1-\frac{a_{Lt}}{a_t}}$ on $t>0$ (see [2, Lemma 3.11 (3.52)]). 

(4) $\ds \frac{t}{\sqrt{T(a_t)}}\sim  Q(a_t)$ and 
$\ds \frac{t\sqrt{T(a_t)}}{a_t}\sim  |Q^\prime(a_t)|$ on $t>0$ (see [2, Lemma 3.4 (3.18) and (3.17)]) .

(5) Assume that $w$ is an Erd\"{o}s-type weight. Then for every $\eta >0$, there exists a constant $C_{\eta} >1$ 
such that
 $$a_{x} \leq C_{\eta} x^{\eta} \ \ (x \geq 1) \leqno(2.2)
 $$
(see [4, Proposition 3 (3.8)]).

\vspace{2ex}

{\bf Lemma 2.2.} ({[2, Theorem 1.9 (a)]}) Let $1 \leq p\le\infty$. Then 
$$\|Pw\|_{L^p(\R)} \le 2\|Pw\|_{L^p([-a_n,a_n])} \leqno(2.3)
$$
for every $n \in \N$ and every $P\in\sa{P}_n$. 

\vspace{2ex}

{\bf Lemma 2.3.}  (1) There exist  constants  $C_1 >1$ and $c_0 >0$ such that if $|x-t| < c_{0}/T(x)$ then  $T(t)/C_{1} \leq T(x) \leq C_{1}T(t)$ holds (cf. [2, Theorem 3.2 (e)] see also [6, Lemma 3.4]).  

(2) There exist a constant $C_{2} >1$ such that for any $n \in \N$, if $|t|, |x| < a_{{2n}}$ and $|x-t| \leq a_{n}/n$ then $T(t)/C_{2} \leq T(x) \leq C_{2}T(t)$ holds (see [6, (4.6)]).

\vspace{2ex}

{\bf Lemma 2.4.} ([5, Theorems 4.1 and 4.2]) Let $m=1,2$ and let $w\in\sa{F}_\lambda(C^{2+m}+)$ with $0 < \lambda < (m+2)/(m+1)$. For every $\alpha\in\R$, 
 we can construct a new weight $w^{*}\in\sa{F}_{\lambda}(C^{1+m}+)$ such that 
$$w^{*}(x) \sim T(x)^{\alpha} w(x) \  \mbox{and} \  T^{*}(x) \sim T(x)\leqno (2.4)$$ on $\R$, and 
$$a_{x/c}\leq a_{x}^{*} \leq a_{cx} \leqno(2.5)$$ 
holds on $\R$ with some constant $c >1$, where $T^{*}$ and $a_{x}^{*}$ are 
corresponding ones defined in 
(1.1) and (1.2) with respect to $w^{*}$ respectively.

\vspace{2ex}

Using the above lemma, we obtain the following assertions. First one is a generalization of (1.7). Second assertion was shown in [5, Corollary 6.2] under some additional assumption. 

\vspace{2ex}

{\bf Lemma 2.5.}  Let $w\in\sa{F}_\lambda(C^{3}+)$ with $0 < \lambda < 3/2$ and let $1 \leq p \leq \infty$.   For $j \in N$, there exists a constant $C_{3} >1$  such that 
for every $n\in \N$ and every $P \in {\sa{P}_{n}}$, we have 
$$ \left\|P^{(j)}\frac{w}{{T^{j/2}}}\right\|_{L^p(\R)}\le C_{3} \kko{\frac{n}{a_n}}^{j} \|Pw\|_{L^p(\R)}
\leqno(2.6) 
$$
and if we further assume that $w \in \sa{F}_\lambda(C^{4}+)$ with $0  < \lambda < 4/3$, then there exists a constant $C_{4} >1$ such that 
$$ \left\|P^{(j)}w\right\|_{L^p(\R)}\le C_{4} \kko{\frac{n}{a_n}}^{j} \| T^{j/2} P w\|_{L^p(\R)}
\leqno(2.7)
$$
also holds.
\vspace{2ex}

{Proof.} For $i =1, \cdots, j$, let $w^{*}_{i} \in \sa{F}_\lambda (C^2 +) $ be a weight obtained in Lemma 2.4  for $\alpha = -(i-1)/2$. Then, since $P^{(j)} \in \sa{P}_{n-j}$,  by (1.7) for $w^{*}_{j}$ and by (2.4) and (2.5), there exists a constant $C>1$ such that
$$\left\|P^{(j)} \frac{w^{*}_{j}}{T^{1/2}}\right\|_{L^p(\R)} \le C \kko{\frac{n-j+1}{a_{(n-j+1)/c}}} \|P^{{(j-1)}}w^{*}_{j}\|_{L^p(\R)}.$$ 
Since $w^{*}_{j}(x) \sim T(x)^{-1/2}w^{*}_{j-1}(x)$, we also see 
$$\left\|P^{(j)} \frac{w^{*}_{j}}{T^{1/2}}\right\|_{L^p(\R)} \leq 
 C\kko{\frac{n-j+1}{a_{(n-j+1)/c}}} \left\|P^{{(j-1)}}\frac{w^{*}_{j-1}}{T^{1/2}}\right\|_{L^p(\R)}.$$
Repeating this process, we have 
\begin{eqnarray*}
& & \left\|P^{(j)}\frac{w}{{T^{j/2}}}\right\|_{L^p(\R)} \leq C\left\|P^{(j)} \frac{w^{*}_{j}}{T^{1/2}}\right\|_{L^p(\R)}\\
& & \ \ \ \le C^{j+1} \kko{\frac{n-j+1}{a_{(n-j+1)/c}}} \cdots \kko{\frac{n}{a_{n/c}}}\left\|Pw\right\|_{L^p(\R)} \\
& & \ \ \ \le C_{3} \kko{\frac{n}{a_{n}}}^{j}\left\|Pw\right\|_{L^p(\R)}, 
\end{eqnarray*}
where we use Lemma 2.1 (1).  

For (2.7), we first remark that if $w \in \sa{F}(C^{3} +)$, then 
$$ \left\|P'w\right\|_{L^p(\R)}\le C_{4} \kko{\frac{n}{a_n}} \| T^{1/2} P w \|_{L^p(\R)}
\leqno(2.8)
$$
holds true
(see [5, Theorem 1.1]). This is the case $j=1$. To show general case $j >1$, we consider a weight $w^{**}_{i} \in \sa{F}_\lambda (C^3 +)$ in 
Lemma 2.4 for $w \in \sa{F}_\lambda (C^4 +)$ with $\alpha = (i-1)/2 \ (i =1, \cdots, j)$. Applying $P^{(j-i)}$ and 
$w^{**}_{i}$ to (2.8) and repeating this process for $i =1, \cdots, j$, we obtain (2.7) as in (2.6).
This completes the proof.

\vspace{2ex}

{\bf Lemma 2.6. } Let $k \in \N \cup \{0\}$ and  $w\in\sa{F}_\lambda(C^{2}+)$ with $0 < \lambda < (k+2)/(k+1)$. Then there exist constants $C_{5} >1$ and $\delta >0$ such that 
$$T(a_{n}) \leq C_{5} n^{2/(2k+3)- \delta } \leqno(2.9)
$$
holds for all $n \in \N$. 

\vspace{2ex}

Proof. We may assume that $w = \exp(-Q)$ is an Erd\"os-type weight. By (2.1), $|Q'(x)|/Q(x)^{\lambda} \leq C$ with some constant $C >1$. 
Hence Lemma 2.1 (4) gives us 
$$\frac{n\sqrt{T(a_n)}}{a_n} \left (\frac{n}{\sqrt{T(a_n)}} \right )^{-\lambda} \leq C, 
$$
that is, $T(a_{n}) \leq C a_{n}^{2/(\lambda +1)} n^{2(\lambda -1)/(\lambda +1)}$. Since $\lambda < (k+2)/(k+1)$, we can choose $\delta >0$ and $\eta >0$ such that $2(\lambda -1)/(\lambda+ 1) + \delta + 2\eta < 2/(2k + 3)$. Hence (2.9) follows from Lemma 2.1 (5). This completes the proof.

\vspace{2ex}

{\bf Lemma 2.7.} Let $w \in \sa{F}_{\lambda}(C^{2} +)$ with $0 < \lambda <2$. Then there exists a constant $C_6 >1$ such that  for every $n \in \N$, if $|t|, |x| < a_{{2n}}$ and if $|t-x| < a_{n}/ (n \sqrt{T(x)})$ then 
$$
w(t)/ C_6 \leq w(x) \leq C_6 w(t)  \leqno(2.10)
$$

\vspace{2ex}

Proof.  By  Lemma 2.3 (2) , we have $T(t)/C_{2} \leq T(x) \leq C_{2}T(t)$, and 
by (1.3) we can write $|t| = a_{s}$. Then 
$a_{s} \leq a_{{2n}}$ implies $s \leq 2n$. Hence (1.4) and Lemma 2.1(1) show $
{s}a_{n}/({n}a_{s}) \leq C_{7}$ with some constant $C_{7} >1$. Since 
$|Q^\prime(t)|\le C s\sqrt{T(a_s)}/{a_s}$ by Lemma 2.1 (4), we have 
\begin{align*}
|Q^\prime(t)||t-x| &\le C\frac{s\sqrt{T(a_s)}}{a_s}\frac{a_n}{n}\frac{1}{\sqrt{T(x)}} \\
&\le C\frac{a_n}{n}\frac{s}{a_{s}} \frac{\sqrt{T(t)}}{\sqrt{T(x)}}
\le C C_{7}\sqrt{C_{2} } .
\end{align*}
 Similarly, we see $|Q'(x)|t-x| \leq C C_{7}$. Hence  if we put $C_6 =  e^{C C_{7}\sqrt{C_{2} }}$, then $|Q^\prime(t)||t-x|\le \log C_6$ and $|Q'(x)||t-x| < \log C_6$ hold true. 
 From mean value theorem for differential calculus, there exists $\theta$ between $x$ and $t$ such that 
$$
\frac{w(x)}{w(t)} = \exp(Q(t)-Q(x)) = \exp(Q^\prime(\theta)(t-x)).$$
Since $Q'$ is increasing, $|Q'(\theta)(x-t)| \leq \max \{|Q'(x)|, |Q'(t)| \} |x-t| \leq \log C_6$, which shows (2.10) immediately. This completes the proof.

\section{Estimates for Christoffel functions}

By definition, the partial sum of Fourier series is given by
\[ s_n(f)(x)=\int_\R K_n(x,t)f(t)w(t)^{2}dt, \leqno{(3.1)} \]
where
\[ K_n(x,t) = \sum_{k=0}^{n-1}p_k(x)p_k(t). \leqno{(3.2)} \]
It is known that by the Cristoffel-Darboux formula
\[ K_n(x,t) = \frac{\gamma_{n-1}}{\gamma_n}\frac{p_n(x)p_{n-1}(t)-p_n(t)p_{n-1}(x)}{x-t} \leqno{(3.3)}  \]
holds, where $\gamma_n$ and $\gamma_{n-1}$ are the leading coefficients of $p_n$ and $p_{n-1}$, respectively.
Then 
\[ a_n \sim \frac{\gamma_{n-1}}{\gamma_n} \leqno(3.4)\]
also holds (see  [2, Lemma 13.9]). 

The Christoffel function $\lambda_{n}(x) = \lambda_n(w,x)$ is defined by
\[ \lambda_n(x):=\frac{1}{K_n(x,x)}=\kko{\sum_{k=0}^{n-1}p_k(x)^{2}}^{-1}. \]
Then 
$$
\lambda_n (x) = \inf_{P\in\sa{P}_{n-1}}\frac{1}{P(x)^{2}} \int_\R |P(t)w(t)|^2dt.
\leqno(3.5)
$$
holds on $\R$.  We use  derivative versions of (3.5). The following equality is also established.

\vspace{2ex}

{\bf Proposition 3.1.} Let $0 \leq j <n$. Then for every $x \in \R$, we have 
$$ \kko{\sum_{k=0}^{n-1}(p_k^{(j)}(x))^{2}}^{-1}  = \inf_{P\in\sa{P}_{n-1}}\frac{1}{(P^{(j)}(x))^{2}} \int_\R |P(t)w(t)|^2dt. \leqno{(3.6)} $$

\vspace{2ex}

Proof. In [3, Theorem 1.3.2], we see
$$
\kko{\sum_{k=0}^{n-1} \Phi(p_{k})^{2}}^{-1}  = \inf_{P\in\sa{P}_{n-1}}\frac{1}{(\Phi( P )^{2}} \int_\R |P(t)w(t)|^2dt
$$
for any linear functional $\Phi$ on polynomials. 
(3.6) follows if we consider $\Phi ( P ) = P^{(j)}(x)$.

\vspace{2ex}

The following estimate plays an important role in our later argument.  We use $C_m  \ ( m= 1,\cdots, 6) $, which are constants in lemmas of the previous section. 

\vspace{2ex}

{\bf Proposition 3.2.} Let $k \geq 2$ be an integer and let $w\in\sa{F}_\lambda(C^{4}+)$ with $0<\lambda<(k+3)/(k+2)$. 
Then there exists a constant $C_{8}>1$ such that for every $ 1 \leq j \leq k$ and every $n \in \N$, 
\[ \frac{w(x)^2}{T(x)^{(2j+1)/2}}\sum_{k=0}^{n-1} (p_k^{(j)}(x))^{2} \le C_{8}\kko{\frac{n}{a_n}}^{2j+1}. \leqno{(3.7)} \]
\vspace{2ex}

Proof. It is enough to show (3.7) for sufficiently large $n$.  By Proposition 3.1, (3.7) follows from $$ \kko{\frac{a_n}{n}}^{2j+1}\frac{w(x)^2}{T(x)^{(2j+1)/2}} \leq C_8 \frac{1}{(P^{(j)}(x))^{2}} \int_\R |P(t)w(t)|^2dt  \leqno(3.8)$$
 for all  $P \in\sa{P}_{n-1}$. Now take  $P\in\sa{P}_{n-1}$ be arbitrarily. By Lemma 2.2, we can choose $\zeta\in \R$ such that  $|\zeta| \leq a_{n-1} $ and satisfies 
$$ \|wP\|_{L^\infty(\R)} \leq 2 |w(\zeta)P(\zeta)|. \leqno(3.9)$$  
Let $0 < c_1 \leq 1$.  Lemma 2.6 gives us  $T(a_n) \leq C_{5}n^{1- \delta'}$ with some $\delta ' >0$, so that  if  $t \in \R$ satisfies $$|t-\zeta| \leq c_1
\frac{a_n}{n} \frac{1}{ \sqrt{T(\zeta)}}, \leqno(3.10) $$then 
\begin{eqnarray*}
& &  |t| \leq |\zeta| + |\zeta -t| \leq  |\zeta|+ c_1 \frac{a_n}{n}\frac{1}{\sqrt{T(\zeta)}} \leq  a_{n-1} + \frac{a_n}{n}  \leq a_n + 
\frac{C_5}{n^{\delta'}} \frac{a_n}{T(a_n)}. 
\end{eqnarray*}
Since  there exists a constant $C >1$ such that  $ a_n + a_n /(C T(a_n)) \leq a_{2n}$ by Lemma 2.1 (3), if we take $n_{0} \in \N$ such that $n_0^{\delta'} > CC_5$, then 
$$|t| \leq a_{2n} \leqno(3.11)$$ for all $n \geq n_0$.  Hence by Lemma 2.7, $ w(t)/C_6 \leq w(\zeta) \leq C_6 w(t)$ holds.  By monotonicity of $w$, 
$w(u)/C_6 \leq w(\zeta) \leq C_6 w(u)$ also holds for every $u$ between $t$ and $\zeta$.  Moreover, since $T$ is quasi-increasing, Lemma 2.3 (2) shows 
 $\sqrt{T(u)}\le C \sqrt{T(\zeta)}$ with some $C >1$.   Then using (2.6) for $p=\infty$ and $j=0$, we have
\begin{align*}
|P(\zeta)|-|P(t)| &\le |P(t)-P(\zeta)| = \ab{\int_\zeta^tP^\prime(u)du} \\
&\le C C_6 \frac{\sqrt{T(\zeta)}}{w(\zeta)}\ab{\int_\zeta^t\frac{1}{\sqrt{T(u)}}w(u)P^\prime(u)du} \\
&\le CC_6 |t-\zeta|\frac{\sqrt{T(\zeta)}}{w(\zeta)}\norm{\frac{w}{\sqrt{T}}P^\prime}{L^\infty(\R)} \\
&\le C C_6 C_3 |t-\zeta|\frac{\sqrt{T(\zeta)}}{w(\zeta)}\frac{n}{a_n}\norm{wP}{L^\infty(\R)} \\
&\le  2 c_1 C C_6 C_3|P(\zeta)|
\end{align*}
by (3.9) and (3.10). Consequently, if we take  $c_1 >0$ so small that $2c_1CC_6C_3 < 1/2$,  we have
$$ |P(t)|\ge\frac{1}{2}|P(\zeta)| \ \ \mbox {if} \ \  |t-\zeta| \leq c_1
\frac{a_n}{n} \frac{1}{ \sqrt{T(\zeta)}}. \leqno(3.12) $$ 
Since $C_2T(t) \geq T(\zeta)$ and $C_6 w(t) \geq w(\zeta)$,  (3.9) and (3.12) show 
\begin{align*}
\int_\R\sqrt{T(t)}|P(t)|^2w(t)^2dt &\ge \frac{\sqrt{T(\zeta)}}{\sqrt{C_2}} \int_{|t-\zeta|\le c_1(a_n/ (n\sqrt{T(\zeta)})}|P(t)|^2w(t)^2dt \\
&\ge \frac{\sqrt{T(\zeta)}}{\sqrt{C_2}} \frac{ |P(\zeta)|^2}{4} \frac{w(\zeta)^2}{C_6^2} c_1 \frac{a_n}{n}\frac{1}{\sqrt{T(\zeta)}}\\
&\geq  \frac{c_1 }{4 \sqrt{C_2}} \frac{1}{C_6^2}\frac{a_n}{n}\frac{\|wP\|_{L^\infty(\R)}^2}{4}\\ & =: \frac{1}{C_{0}}\frac{a_n}{n}\|wP\|_{L^\infty(\R)}^2.
\end{align*}
We note that in the above argument we only use the fact that $w \in \sa{F}_\lambda (C^3 +)$. If  $w \in \sa{F}_\lambda(C^4 +)$,  we  can construct $w^* 
\in \sa{F}_\lambda(C^3 + )$ such that $w^*(x) \sim T(x)^{-1/4} w(x)$ by Lemma 2.4.  Then it follows from (2.6) for $p=\infty$  that for every $x \in R$, 
\begin{align*}
\int_\R\sqrt{T^*(t)}|P(t)|^2w^*(t)^2dt &\ge \frac{1}{C_{0}} \frac{a_n^*}{n}\|w^*P\|_{L^\infty(\R)}^2 \\
&\ge \frac{1}{C_{0}C_3} \frac{a_{n}^{*}}{n} \kko{\frac{a_{n-1}^*}{n-1}}^{2j}\norm{\frac{w^*}{(T^{*})^{j/2}}P^{(j)}}{L^\infty(\R)}^2 \\
&\ge \frac{1} {C_{0} C_3} \kko{\frac{a_n^*}{n}}^{2j+1} \frac{w^*(x)^2}{T^*(x)^{j}}|P^{(j)}(x)|^2, 
\end{align*}
and hence by (2.4) and (2.5) we see 
\begin{align*}
\int_\R|P(t)|^2w^2(t)dt & \geq \frac{1}{C}  \int_\R\sqrt{T^\ast(t)}|P(t)|^2{w^\ast(t)}^2dt \\
&\ge \frac{1}{CC_{0} C_3} \kko{\frac{a_n^*}{n}}^{2j +1 }\frac{w^*(x)^2}{T^*(x)^{j}}|P^{(j)}(x)|^2\\
& \geq  \frac{1}{C} \kko{\frac{a_{n/c}}{n}}^{2j +1 }\frac{w(x)^2}{T(x)^{(2j +1)/2}}|P^{(j)}(x)|^2.
\end{align*}
This together with Lemma 2.1 (1) shows (3.8) and the proof is completed. 

\section{Proof of Theorem 1.1}

In the remaining sections, we again use $C_m  \ (m= 1,\cdots, 6) $ without notice, which are constants in Lemmas of the previous section. 

 Let  $1 \leq p \leq \infty$, $k \geq 2, w \in \sa{F}_{\lambda}(C^{4} +)$ with 
$0 < \lambda <(k+3)/(k+2)$ and let $1 \leq j \leq k$.  Due to Lemma 2.4, there is $w^{*} \in \sa{F}_{\lambda}(C^{3} + )$ 
such that $w^{*}(x) \sim T(x)^{-(2j+1)/4}w(x)$. Take $fw \in L^{p}(\R)$ arbitrarily. Since $v_n^{(j)}(f)\in\sa{P}_{2n-j}$, applying $w^{*}$ to 
 (2.7), we have
\begin{eqnarray*}
& & \left\|v_n^{(j)}(f)\frac{w}{T^{(2j+1)/4}}\right\|_{L^p(\R)} \le C\|v_n^{(j)}(f)w^\ast \|_{L^p(\R)} \\
& & \ \ \ \le C\kko{\frac{2n -j}{a^{*}_{2n-j}}}^j\|(T^{*})^{j/2}v_n(f)w^\ast\|_{L^p(\R)} \\
& & \ \ \ \le C\kko{\frac{n}{a_{2n/c}}}^j\left\|v_n(f)\frac{w}{T^{1/4}}\right\|_{L^p(\R)}\\
& & \ \ \ \le C\kko{\frac{n}{a_n}}^j\|fw\|_{L^p(\R)}.
\end{eqnarray*}
Here we use Lemma 2.1 (1),  (2.4) and (2.5). The last inequality follows from (1.5). This completes the proof of Theorem 1.1.

\vspace{2ex}

By a similar argument as above, we also have
\[ \|v_n^{{(j)}}(f) w\|_{L^p(\R)} \le C\kko{\frac{n}{a_n}}^jT(a_n)^{(2j+1)/4}\|fw\|_{L^p(\R)}.  \leqno{(4.1)} \]

In fact, take $w^{*} \in \sa{F}_{\lambda}(C^{3} +)$ such that $w^{*}(x) \sim T^{j/2}(x)w(x)$. Then by (2.7) for $w$ and by
Lemma 2.4 for $w^{*}$, we have 
\begin{eqnarray*}
& & \|v_n^{(j)}(f)w\|_{L^p(\R)} \le C\kko{\frac{n}{a_n}}^j\|T^{j/2}v_n(f)w \|_{L^p(\R)} \\
& & \ \ \ \le C\kko{\frac{n}{a_n}}^j\|v_n(f)w^\ast\|_{L^p([-a_{2n}^\ast,a_{2n}^\ast])} \\
& & \ \ \ \le C\kko{\frac{n}{a_n}}^j\left\|v_n(f)T^{(2j+1)/4}\frac{w}{T^{1/4}}\right\|_{L^p([-a_{2cn},a_{2cn}])} \\
& & \ \ \ \le C\kko{\frac{n}{a_n}}^jT(a_n)^{(2j+1)/4}\left\|v_n(f)\frac{w}{T^{1/4}}\right\|_{L^p([-a_{2cn},a_{2cn}])}\\
& & \ \ \ \le C\kko{\frac{n}{a_n}}^jT(a_n)^{(2j+1)/4}\|fw\|_{L^p(\R)}. 
\end{eqnarray*}
Note that by Lemma 2.1 (2), $T(x) \leq CT(a_{2cn}) \leq CT(a_{n}) $ holds for all $x \in [-a_{2cn},a_{2cn}]$, because $T$ is quasi-increasing. 

\section{Proof of Theorem 1.2}

 Let $k \geq 2, w \in \sa{F}_{\lambda}(C^{4} +)$ with 
$0 < \lambda <(k+3)/(k+2)$ and let $1 \leq j \leq k$. We first show (1.12) for the case $p=\infty$. 
Suppose that $T^{(2j+1)/4}fw \in L^\infty(\R)$. Since $v_n^{(j)}(f) \in \sa{P}_{2n}$, by Lemma 2.2, it is sufficient to show 
$$  |v_n^{{(j)}}(f)(x)  w (x) | \le C\kko{\frac{n}{a_n}}^j\|T^{(2j+1)/4}fw\|_{L^\infty(\R)}  \leqno{(5.1)}$$
for every $|x| \leq a_{2n}$. 
Now we set 
$$
A_n  :=   \{t \in \R; |t-x| < \frac{a_{2n}}{2n} \}, \ \ 
B_n  :=   \{t \in \R;  \frac{a_{2n}}{2n} \leq  |t-x| < \frac{c_0}{T(x)} \} $$
and $ C_n  :=  \R \setminus (A_n  \cup B_n)$, where $c_0 >0$ is a constant in Lemma 2.3 (1).  Then as in the proof of (3.11), there exists $n_0 \in N$ 
such that if $n \geq n_0$ and $t \in A_n$, then $|t| \leq a_{4n}$ holds. Hence Lemma 2.3 (2) implies 
$$
T(t)/C_2\leq T(x) \leq C_2 T(t) \ \ (t \in A_n). \leqno(5.2)
$$ 
Since $T$ is bounded on $[-a_{4n_0}, a_{4n_0}]$, we may assume that (5.2) holds for all $n \in \N$.  Also by Lemma 2.3 (1), 
$$
T(t)/C_1 \leq T(x) \leq C_1 T(t) \ \ (t \in B_n) \leqno(5.3)
$$
holds true.  Let $g(t) := f(t) \chi_{A_n}(t)$, where $\chi_A$ is the characteristic function of a set $A$ and put $h(t) = f(t) - g(t)$.  Since 
\[ \int_\R\kko{\sum_{k=0}^{m-1}p_k^{(j)}(x)p_k(t)}^2w(t)^2dt=\sum_{k=0}^{m-1} (p_k^{(j)}(x))^2, \]
(3.2), (5.2) and the Schwarz inequality show that 
\begin{align*}
&|s_m^{(j)}(g)(x)w(x)| \\
&\le w(x)\int_\R \ab{g(t)\sum_{k=0}^{m-1}p_k^{(j)}(x)p_k(t)w(t)^2}dt \\
&\le \kko{\sum_{k=0}^{m-1}(p_k^{(j)}(x))^2 w(x)^2}^{1/2}\kko{\int_{A_{n}}|f(t)w(t)|^2dt}^{1/2} \\
&\le C_2^{(2j+1)/4}  \left ( \sum_{k=0}^{m-1}\frac{w(x)^2}{T(x) ^{(2j +1)/2}}({p_k^{(j)}}(x))^2 \right )^{1/2} \kko{\int_{A_{n}}|T(t)^{(2j+1)/4}f(t)w(t)|^2dt}^{1/2} \\
&\le C\left ( \sum_{k=0}^{m-1}\frac{w(x)^2}{T(x) ^{(2j +1)/2}}({p_k^{(j)}}(x))^2 \right )^{1/2}
\| T^{(2j+1)/4}f w \|_{L^\infty(\R)}  \left (  \frac{a_{2n}}{2n} \right )^{1/2} .
\end{align*}
Since $v_n^{(j)}(g) (x) =  (1/n)
\sum^{2n}_{m = n+1} s_m^{(j)}(g)(x)$, Proposition 3.2 gives us 
\[ |v_n^{(j)}(g)(x)w(x)|\le C\kko{\frac{n}{a_n}}^j\|T^{(2j+1)/4}fw\|_{L^\infty(\R)} \leqno(5.4) \]
for all $x \in R$ with $|x| \leq a_{2n}$.

To estimate $v_{n}^{(j)}(h)$, we use (3.3).  For $i = 0,1, \cdots, j$, we put
\begin{align*}
&v_{n,i}(h)(x)\\
:&= \frac{1}{n}\sum_{m=n+1}^{2n}\frac{\gamma_{m-1}}{\gamma_m}\int_\R h(t)\frac{p_m^{(j-i)}(x)p_{m-1}(t)-p_{m-1}^{(j-i)}(x)p_m(t)}{(x-t)^{i+1}}w(t)^{2}dt \\
&= \frac{1}{n}\sum_{m=n+1}^{2n}\frac{\gamma_{m-1}}{\gamma_m}(b_{m-1}(h_i)p_m^{(j-i)}(x)-b_m(h_i)p_{m-1}^{(j-i)}(x)),
\end{align*}
where
$$
h_i(t):=\frac{h(t)}{(x-t)^{i+1}} \ \ \mbox{and} \ \  b_k(h_i):=\int_\R h_i(t)p_k(t)w(t)^{2}dt \ \ (k \in \N \cup \{0\}).
$$
Then
$$ v_n^{(j)} (h)(x)=\sum_{i=0}^j (-1)^i\binom{j}{i} v_{n,i}(h)(x). \leqno(5.5) $$
By (3.4), the Schwarz inequality  and Proposition 3.2, we have 
\begin{align*}
&|v_{n,i}(h)(x)w(x)| \\
&\le \frac{1}{n}\sum_{m=0}^{2n}\ab{\frac{\gamma_{m-1}}{\gamma_m}2p_m^{(j-i)} (x)b_m(h_i)w(x)} \\
&\le C\frac{a_n}{n}\kko{w(x)^{2}\sum_{m=0}^{2n}({{p_m^{(j-i)}}(x))^{2}}}^{1/2}\kko{\sum_{m=0}^{2n}|b_m(h_i)|^2}^{1/2} \\
&\le C\frac{a_n}{n}\kko{\frac{w(x)^{2}}{T(x)^{(2(j-i)+1)/2}}\sum_{m=0}^{2n}({{p_m^{(j-i)}}(x)})^{2}}^{1/2}\kko{T(x)^{(2(j-i)+1)/2}\sum_{m=0}^{2n}|b_m(h_i)|^2}^{1/2} \\
&\le C\kko{\frac{n}{a_n}}^{(2(j-i)-1)/2}\kko{T(x)^{(2(j-i)+1)/2}\sum_{m=0}^{2n}|b_m(h_i)|^2}^{1/2}. 
\end{align*}
The Bessel inequality implies that  
\begin{align*}
\sum_{m=0}^{2n}|b_m(h_i)|^2 \le \int_\R\ab{\frac{h(t)}{(x-t)^{i+1}}}^2w(t)^2dt = \int_{B_{n} \cup C_{n}}\frac{|f(t)w(t)|^2}{(x-t)^{2(i+1)}}dt
\end{align*}
and hence, by (5.3), we have
\begin{align*}
& T(x)^{(2(j-i)+1)/2}\int_{B_{n}}\frac{|f(t)w(t)|^2}{(x-t)^{2(i+1)}}dt \\
&\le C_{1}^{(2(j-i)+1)/2} \int_{B_{n}}\frac{|T(t)^{(2(j-i)+1)/4}f(t)w(t)|^2}{(x-t)^{2(i+1)}}dt \\
&\le C\|T^{(2(j-i)+1)/4}fw\|_{L^\infty(\R)}^2\int_{|x-t| > \frac{a_{{2n}}}{2n}}\frac{1}{(x-t)^{2(i+1)}}dt \\
&\le C\|T^{(2(j-i)+1)/4}fw\|_{L^\infty(\R)}^2\kko{\frac{n}{a_{n}}}^{2i+1}\\
&\le C\|T^{(2j+1)/4}fw\|_{L^\infty(\R)}^2\kko{\frac{n}{a_{n}}}^{2i+1},
\end{align*}
because $T \geq 1$. On the other hand, if $|x| \leq a_{{2n}}$ then $T(x) \leq C T(a_{n})$,  so that  
\begin{align*}
& T(x)^{(2(j-i)+1)/2}\int_{C_{n}}\frac{|f(t)w(t)|^2}{(x-t)^{2(i+1)}}dt \\
&\le C\|fw\|_{L^\infty(\R)}^2T(x)^{(2(j-i)+1)/2}\int_{\frac{c_{0}}{T(x)}\le|x-t|}\frac{1}{(x-t)^{2(i+1)}}dt \\
&\le C\|fw\|_{L^\infty(\R)}^2T(x)^{(2(i+j)+3)/2} \\
& \le C \|T^{(2j+1)/4}fw\|_{L^{\infty}(\R)}^{2} T(a_{n})^{(2(i+k) +3)/2}.
\end{align*}
Moreover
$$
T(a_{n})^{(2(i+k) +3)/2} \leq C \left ( \frac{n}{a_{n}} \right )^{2i +1} \leqno(5.6)
$$
holds. In fact, to show this we may assume that $w$ is an Erd{\"{o}}s-type weight. Then by Lemma 2.1 (5) and Lemma 2.6, 
we have 
\begin{align*}
T^{(2k+3)/2}(a_{n}) & \le C n^{(2/(2k+3)-\delta)((2k+3)/2)} \le Cn^{1-\delta'}\le C\kko{\frac{n}{a_n}}.
\end{align*}
Similarly 
\begin{align*}
T^{(2(i+k)+3)/2}(x) &\le CT(a_n)^{(4k+3)/2} \le C n^{(2/(2k+3)-\delta)((4k+3)/2)} \\
&\le Cn^{2-\delta''}\le C\kko{\frac{n}{a_n}}^2\le C\kko{\frac{n}{a_n}}^{2i+1}
\end{align*}
holds for $i \geq 1$. 
Combining the above estimates, we thus have 
\begin{align*}
|v_{{n,i}}(h)(x)w(x)| & \leq C\kko{\frac{n}{a_n}}^{(2(j-i)-1)/2}\kko{T(x)^{(2(j-i)+1)/2}\sum_{m=0}^{2n}|b_m(h_i)|^2}^{1/2}
\\ & 
\leq C\kko{\frac{n}{a_n}}^{(2(j-i)-1)/2} \|T^{(2j+1)/4}fw\|_{L^{\infty}(\R)} \left ( \frac{n}{a_{n}} \right )^{(2i +1)/2}\\
\\ & \leq C\kko{\frac{n}{a_n}}^{j} \|T^{(2j+1)/4}fw\|_{L^{\infty}(\R)}. 
\end{align*} 
It follows from (5.5) that 
$$|v_n^{{(j)}}(h)(x)  w (x) | \le C\kko{\frac{n}{a_n}}^j\|T^{(2j+1)/4}fw\|_{L^\infty(\R)}. $$ This together with (5.4) shows  (5.1). 

We will prove (1.12) for $p =2$ in the next  section.  
Then using the Riesz-Thorin interpolation theorem for an operator
\[ F~:~f~\mapsto~wv_n^{(j)}\kko{\frac{f}{w}}, \]
we obtain (1.12) for all $2\le p\le\infty$. This completes the proof of Theorem 1.2.

\section{Proof of Theorem 1.3}

Let $1 \leq p \leq2$ and $T^{(2j+1)}fw \in L^{2}(\R)$. We use the same notations as in the previous section.  
Then as in the estimate of $s_{m}^{(j)}(g)$ in the previous section, we have 
$$|s_m^{(j)}(g)(x)w(x)| \leq C \left (  \frac{n}{a_{n}} \right )^{(2j+1)/2} \kko{\int_{A_{n}}|
T(t)^{(2j+1)/4}f(t)w(t)|^2dt}^{1/2} \leqno(6.1)$$
for $|x| \leq a_{2n}$. Hence Lemma 2.2 and the H\"older inequality imply
\begin{align*}
& \int_\R|s_m^{(j)}(g)(x)w(x)|^pdx  \leq 2^{p} \int_{|x|\leq a_{2n}}|s_m^{(j)}(g)(x)w(x)|^pdx \\
&\le C\int_{|x| \leq a_{2n}} \kko{\frac{n}{a_n}}^{p(2j+1)/2}\kko{\int_{A_{n}}|T(t)^{(2j+1)/4}f(t)w(t)|^2dt}^{p/2}dx \\
& \leq C\kko{\frac{n}{a_n}}^{p(2j+1)/2}\int_{|x| \leq a_{2n}} \kko{\int_{|u|\le\frac{a_{2n}}{2n}} |T(x-u)^{(2j+1)/4}f(x-u)w(x-u)|^2du}^{p/2}dx \\
&\le C\kko{\frac{n}{a_n}}^{p(2j+1)/2}a_n^{(2-p)/2}\\
&~~~~\times\kkko{\int_{|x| \leq a_{2n}}\kko{\int_{|u|\le\frac{a_n}{n}}|T(x-u)^{(2j+1)/4}f(x-u)w(x-u)|^2du}dx}^{p/2} \\
&\leq C\kko{\frac{n}{a_n}}^{p(2j+1)/2}a_n^{(2-p)/2}\|T^{(2j+1)/4}fw\|_{L^2(\R)}^p\kko{\int_{|u|\le\frac{a_n}{n}}du}^{p/2} \\
&\leq C\kko{\frac{n}{a_n}}^{pj}a_n^{(2-p)/2}\|T^{(2j+1)/4}fw\|_{L^2(\R)}^p,
\end{align*}
so that we have
$$ \|v_n^{(j)}(g)w\|_{L^p(\R)}\le C\kko{\frac{n}{a_n}}^{j}a_n^{(2-p)/2p}\|T^{(2j+1)/4}fw\|_{L^2(\R)}.\leqno(6.2)$$
Next we estimate $v_{n,i}(h)$. Similarly as above, we have
\begin{align*}
&\int_\R|v_{n,i}(h)(x)w(x)|^pdx \leq 2 \int_{|x| \leq a_{2n}} |v_{n,i}(h)(x)w(x)|^{p}dx\\
&\le C\kko{\frac{n}{a_n}}^{p(2(j-i)-1)/2}a_n^{(2-p)/2}\\
&~~~~\times\kkko{\int_{|x|\leq a_{2n}}\kko{\int_{B_{n} \cup C_{n}}\frac{|T(t)^{(2(j-i)+1)/4}f(t)w(t)|^2}{(x-t)^{2(i+1)}}dt}dx}^{p/2}.
\end{align*}
Also as in the argument of previous section, 
\begin{align*}
&\int_{|x| \leq a_{2n}} \kko{\int_{B_{n}}\frac{|T^{(2(j-i)+1)/4}(t)f(t)w(t)|^2}{(x-t)^{2(i+1)}}dt}dx \\
&\le \int_\R\kko{\int_{\frac{a_n}{n}\le|u|}\frac{|T^{(2(j-i)+1)/4}(x-u)f(x-u)w(x-u)|^2}{u^{2(i+1)}}du}dx \\
&\le C\kko{\frac{n}{a_n}}^{2i+1}\|T^{(2(j-i)+1)/4}fw\|_{L^2(\R)}^2 \leq   C\kko{\frac{n}{a_n}}^{2i+1}\|T^{(2j+1)/4}fw\|_{L^2(\R)}^2
\end{align*}
On the other hand, by (5.6) we have
\begin{align*}
&\int_{|x| \leq a_{2n}} \kko{T(x) ^{(2(j-i)+1)/2}\int_{C_{n}}\frac{|f(t)w(t)|^2}{(x-t)^{2(i+1)}}dt}dx \\
&\le C T(a_{2n})^{(2(j-i)+1)/2}\int_\R\kko{\int_{\frac{c_{0}}{T(a_{2n})}\le|u|}\frac{|f(x-u)w(x-u)|^2}{u^{2(i+1)}}du}dx \\
&\le C\|fw\|_{L^2(\R)}^2T(a_{2n})^{(2(j-i)+1)/2}\int_{\frac{c_{0}}{T(a_{2n})}\le|u|}\frac{1}{u^{2(i+1)}}du \\
&\le CT(a_{2n})^{(2j+2i+3)/2}\|fw\|_{L^2(\R)}^2 \\
&\le C\kko{\frac{n}{a_n}}^{2i+1}\|T^{(2j+1)/4}fw\|_{L^2(\R)}^2. 
\end{align*}
Consequently  we have
$$ \|v_{n,i}(h)w\|_{L^p(\R)}\le C\kko{\frac{n}{a_n}}^{j}a_n^{(2-p)/2p}\|T^{(2j+1)/4}fw\|_{L^2(\R)} \leqno(6.3)$$
for $0 \leq i \leq j$, so that 
\[ \|v_n^{(j)}(h)w\|_{L^p(\R)}\le C\kko{\frac{n}{a_n}}^j a_n^{(2-p)/2p}\|T^{(2j+1)/4}fw\|_{L^2(\R)}\]
follows. This together with (6.2) shows (1.13). This completes the proof of Theorem 1.3.

\vspace{2ex}

Under the same assumptions in Theorem 1.3, the following estimate is also established. Let $\beta >1$ and $1 \leq p \leq 2$. Then
$$ \|v_n^{(j)}(f) \frac{w}{(1+|x|)^{(2-p)\beta/(2p)}}\|_{L^p(\R)}\le C\kko{\frac{n}{a_n}}^{j}\|T^{(2j+1)/4}fw\|_{L^2(\R)}\leqno(6.4)$$
holds for every $T^{(2j+1)/4}fw \in L^{2}(\R)$ and every $n \in \N$. 

In fact, in the proof of Theorem 1.3, we used
\begin{eqnarray*}
& & \int_{|x| \leq a_{2n}} \kko{\int_{|x-t|\le\frac{a_{2n}}{2n}} |T(t)^{(2j+1)/4}f(t)w(t)|^2dt}^{p/2}dx \\
& & \ \ \ \leq a_n^{(2-p)/2}
\kkko{\int_{|x| \leq a_{2n}}\kko{\int_{|x-t|\le\frac{a_n}{n}}|T(t)^{(2j+1)/4}f(t)w(t)|^2du}dx}^{p/2}, 
\end{eqnarray*}
which follows from the H\"older inequality. Instead of this, we use
\begin{eqnarray*}
& & \int_{\R}\frac{1}{(1+|x|)^{(2-p)\beta/2}} \kko{\int_{|x-t|\le\frac{a_{2n}}{2n}} |T(t)^{(2j+1)/4}f(t)w(t)|^2dt}^{p/2}dx \\
& & \ \ \leq \left ( \int_{\R}\frac{1}{(1+|x|)^{\beta}} dx \right )^{(2-p)/2} 
\kkko{\int_{\R}\kko{\int_{|x-t|\le\frac{a_n}{n}}|T(t)^{(2j+1)/4}f(t)w(t)|^2dt}dx}^{p/2}.
\end{eqnarray*}
Then as in (6.2), we obtain 
$$ \|v_n^{(j)}(g) \frac{w}{(1+|x|)^{(2-p)\beta/(2p)}}\|_{L^p(\R)}\le C\kko{\frac{n}{a_n}}^{j}\|T^{(2j+1)/4}fw\|_{L^2(\R)}.$$
For the estimate of $v_{n,i}(h)$, we take $w^{*} \in \sa{F}_{\lambda}(C^{3} +)$ such that 
$w^{*}(x) \sim w(x)/(1+ |x|)^{(2-p)\beta/(2p)}$ (see [5, Theorem 4.2]). Then by Lemma 2.2, 
$$
\int_{R} \left | v_{n,i}(h) \frac{w(x)}{(1+|x|)^{(2-p)\beta/(2p)}} \right|^{p} dx \leq 2^{p}
\int_{|x| \leq a^{*}_{2n}} \left | v_{n,i}(h) \frac{w(x)}{(1+|x|)^{(2-p)\beta/(2p)}} \right|^{p} dx. 
$$
By an estimate similar to (6.3), we obtain 
$$ \|v_{n,i}(h) \frac{w}{(1+|x|)^{(2-p)\beta/(2p)}}\|_{L^p(\R)}\le C\kko{\frac{n}{a_n}}^{j}\|T^{(2j+1)/4}fw\|_{L^2(\R)},$$
which shows (6.4).


\mbox{}\\
{\footnotesize{
\qw Kentaro Itoh\\
Department of Mathematics\\
Meijo University\\
Tenpaku-ku, Nagoya 468-8502\\
Aichi Japan\\
{\texttt{133451501@ccalumni.meijo-u.ac.jp}}
\vspace{1ex}
\qw Ryozi Sakai\\
Department of Mathematics\\
Meijo University\\
Tenpaku-ku, Nagoya 468-8502\\
Aichi Japan\\
{\texttt{ryozi@crest.ocn.ne.jp}}
\vspace{1ex}
\qw Noriaki Suzuki\\
Department of Mathematics\\
Meijo University\\
Tenpaku-ku, Nagoya 468-8502\\
Aichi Japan\\
{\texttt{suzukin@meijo-u.ac.jp}}
}}

\begin{thebibliography}{99}
\bibitem{Ein01} H. S. Jung and R. Sakai, Specific examples of exponential weights,
Commun. Korean Math. Soc. 24. No.2 (2009), pp.303-319.
\bibitem{Ein02} A. L. Levin and D. S. Lubinsky, Orthogonal Polynomials for Exponential
Weights, Springer, New York, 2001.
\bibitem{Ein03} H. N. Mhaskar, Introduction to the Theory of Weighted Polynomial
Approximation, World Scientific, Singapore, 1996.
\bibitem{Ein04} R. Sakai and N. Suzuki, Favard-type inequalities for 
exponential weights, Pioneer J. of Math. Vol 3, No.1 (2011), pp.1-16.
\bibitem{Ein05} R. Sakai and N. Suzuki, Mollification of exponential weights and its
application to the Markov-Bernstein inequality, Pioneer J. of Math., Vol.7,
No.1 (2013), pp.83-101.
\bibitem{Ein06} K. Itoh, R. Sakai and N. Suzuki, The de la Vall\'{e}e Poussin mean and polynomial approximation for exponential weight,  arXiv:1311.3337 [math.CA].
\end{thebibliography}
\end{document}